\documentclass{amsart}
\usepackage[all]{xy}
\usepackage{latexsym}
\usepackage{amssymb}
\usepackage{amsmath}
\usepackage{amsthm}
\usepackage{amscd}
\usepackage{fancyhdr}
\usepackage[dvips]{graphicx}
\usepackage{fullpage}
\usepackage{mathrsfs}
\input cyracc.def
\xyoption{arc}
 \newfam\cyrfam

  \font\tencyr=wncyr10

  \font\sevencyr=wncyr7

  \font\fivecyr=wncyr5

  \textfont\cyrfam=\tencyr \scriptfont\cyrfam=\sevencyr

    \scriptscriptfont\cyrfam=\fivecyr

  \newfam\cyifam

  \font\tencyi=wncyi10

  \font\sevencyi=wncyi7

  \font\fivecyi=wncyi5

  \textfont\cyifam=\tencyi \scriptfont\cyifam=\sevencyi

    \scriptscriptfont\cyifam=\fivecyi

 \newcommand{\lon}{\longrightarrow}
 \newcommand{\bu}{\bullet}
 
 \newcommand{\rar}{\rightarrow}

\newcommand{\p}{{\partial}}
\newcommand{\Id}{{\mathrm I\mathrm d}}

 \newcommand{\Z}{{\mathbb Z}}
 \newcommand{\bS}{{\mathbb S}}

 \newcommand{\K}{{\mathbb K}}

\newcommand{\alg}[1]{\mathfrak{#1}}
 \newcommand{\ot}{\otimes}

\newcommand{\sC}{{\mathsf C}}

\newcommand{\sE}{{\mathsf E}}

\newcommand{\sG}{{\mathsf G}}
\newcommand{\sP}{{\mathsf P}}

\newcommand{\sa}{{\mathsf a}}
\newcommand{\ssf}{{\mathsf f}}
\newcommand{\sr}{{\mathsf r}}

\newcommand{\Ber}{{\mathit B}{\mathit e} {\mathit r}}

 \newcommand{\Beq}{\begin{equation}}
 \newcommand{\Eeq}{\end{equation}}
 \newcommand{\Beqr}{\begin{eqnarray}}
 \newcommand{\Eeqr}{\end{eqnarray}}
 \newcommand{\Beqrn}{\begin{eqnarray*}}
 \newcommand{\Eeqrn}{\end{eqnarray*}}
 \newcommand{\Ba}{\begin{array}}
 \newcommand{\Ea}{\end{array}}
 \newcommand{\Bi}{\begin{itemize}}
 \newcommand{\Ei}{\end{itemize}}
 \newcommand{\Bc}{\begin{center}}
 \newcommand{\Ec}{\end{center}}
 \newcommand{\fg}{{\mathfrak g}}

 \newcommand{\f}{{\mathcal O}}

 \newcommand{\cC}{{\mathcal C}}
 
 \newcommand{\cE}{{\mathcal E}}

 \newcommand{\cM}{{\mathcal M}}
 
 \newcommand{\cP}{{\mathcal P}}
 \newcommand{\cQ}{{\mathcal Q}}

 \newcommand{\ga}{\gamma}
 
 \newcommand{\Ga}{\Gamma}

 \newcommand{\om}{\omega}

 \newcommand{\Hom}{{\mathrm H\mathrm o\mathrm m}}
 
 \newcommand{\sip}{\smallskip}
 \newcommand{\bip}{\bigskip}
 
\theoremstyle{plain}
\swapnumbers
\newtheorem{theorem}{Theorem}[subsection]

\newtheorem{proposition}[theorem]{Proposition}

\newtheorem{prop-def}[theorem]{Proposition-definition}

\newtheorem{f-theorem}{Formality Theorem}[section]
\newtheorem{main-theorem}{Main~Theorem}[section]
\newtheorem{section-theorem}{Theorem}[section]
\newtheorem{section-corollary}{Corollary}[section]

\theoremstyle{definition}

\newtheorem{remark}[theorem]{Remark}

\newtheorem{fact-me}{Fact \cite{Me1}}[subsection]

\newcommand{\grt}{\alg{grt}_1}

\newcommand{\GRT}{GRT_1}

\newcommand{\hu}{u}

\begin{document}

 \sloppy

 \newenvironment{proo}{\begin{trivlist} \item{\sc {Proof.}}}
  {\hfill $\square$ \end{trivlist}}

\long\def\symbolfootnote[#1]#2{\begingroup%
\def\thefootnote{\fnsymbol{footnote}}\footnote[#1]{#2}\endgroup}

 \title{Grothendieck-Teichm\"uller  and Batalin-Vilkovisky}
  %\title{GRT and BV}

 \author{ Sergei\ Merkulov and Thomas Willwacher}
\address{Sergei~Merkulov: Department of Mathematics, Stockholm University, 10691 Stockholm, Sweden, and Mathematics Research Unit, University of Luxembourg, Grand Duchy of Luxembourg (Current address)}
\email{sergei.merkulov@uni.lu}
%\address{Thomas~Willwacher: Department of Mathematics, Harvard University}
\address{Institute of Mathematics\\ University of Zurich\\ Winterthurerstrasse 190 \\ 8057 Zurich, Switzerland}
\email{thomas.willwacher@math.uzh.ch}


 \begin{abstract} It is proven that, for  any affine supermanifold $M$  equipped with a constant odd symplectic structure, there is a universal
  action (up to homotopy) of the  Grothendieck-Teichm\"uller Lie algebra $\grt$ on the set of quantum  BV structures (i. e.\ solutions of the quantum master equation) on $M$.

%\sip
%\noindent {\sc Mathematics Subject Classifications} (2000). %53D55, 16E40, 18G55, 58A50.

%\noindent {\sc Key words}. Kontsevich graph complex, BV formalism.
\end{abstract}
 \maketitle
\markboth{S.\ Merkulov and T.\ Willwacher}{GRT and BV}

{\large
\section{\bf Introduction}
}

Let $M$ be a finite dimensional affine $\Z$-graded manifold $M$  over a field $\K$ equipped with a constant degree 1 symplectic structure
$\om$. In particular, the ring of functions $\f_M$ is a Batalin-Vilkovisky algebra, with Batalin-Vilkovisky operator $\Delta$ and bracket $\{\ ,\ \}$. A degree 2 function $S\in \f_M[[\hu]]$ is a solution the quantum master equation on $M$ if\footnote{See \cite{Sc} for an introduction into the geometry of the BV formalism.}
$$
\hu \Delta S + \frac{1}{2}\{S,S\}=0,
$$
where $u$ is a formal variable of degree 2. In other words $S$ is a Maurer-Cartan element in the differential graded (dg) Lie algebra $\left(\f_M[[\hu]][1], \hu\Delta, \{\ ,\ \}\right)$.

\sip

The Grothendieck-Teichm\"uller group $\GRT$ is a pro-unipotent group  introduced
by Drinfeld in [Dr]; we denote its Lie algebra by $\grt$.
In this paper we show the following result.

 \sip

\subsection{ Main Theorem} {\em There is an $L_\infty$ action of the Lie algebra $\grt$ on the differential graded Lie algebra $\left(\f_M[[\hu]][1], \hu\Delta, \{\ ,\ \}\right)$ by $L_\infty$ derivations. In particular, it follows that there is an action of $\GRT$ on the set of gauge equivalence classes of formal solutions of the quantum master equation, i.~e., on gauge equivalence classes of Maurer-Cartan elements in the differential graded Lie algebra $\left(\hbar\f_M[[\hu]][[\hbar]][1], \hu\Delta, \{\ ,\ \}\right)$, where $\hbar$ is a formal deformation parameter of degree 0. }

 \sip

 Our main technical tool is a  version of the Kontsevich graph complex, $(\sG\sC_2[[\hu]], d_\hu)$ which
 controls universal deformations of $\left(\f_M[[\hu]][1], \hu\Delta, \{\ ,\ \}\right)$  in the category of $L_\infty$ algebras.
 Using the main result of \cite{Wi} we show in Sect.\ 2 that
 $$
 H^0(\sG\sC_2[[\hu]], d_\hu)\simeq\grt
 $$
and then use this isomorphism in Sect.\ 3 to prove the Main Theorem.

%  In Sect.\ 3 we explain how to use this isomorphism of Lie algebras to define a universal homotopy action of $\grt$  on the set of quantum BV structures on any affine odd symplectic manifold $M$. More precisely, we prove in Sect.\ 3 the following
%  \bip
%
%
%  \subsection{ Main Theorem} {\em There is a universal morphism of groups,
% $$
% GRT_1 \lon Ho AutL_\infty(M)
% $$
% where $Ho Aut L_\infty(M)$ is the group of homotopy classes of
% $L_\infty$ automorphisms of the dg Lie algebra $(\f_M[[\hu]], \hu\Delta, \{\ ,\ \})$. }
%
% \bip
%
% The precise definition of  the group $Ho Aut L_\infty(M)$ is given in Sect. \ref{} below with the help
% of the Deligne-Getzler 2-groupoid construction \cite{De,Ge}. The word {\em universal}\, means that the above morphism makes sense for {\em any}\, formal finite-dimensional manifold equipped with a degree 1 symplectic form. We would like to emphasize that the finite-dimensionality assumption is important, otherwise the formulae lying behind the above morphism of groups might give divergent results.

\subsection{ Some notation} In this paper $\mathbb K$ denotes a field of characteristic $0$. If $V=\oplus_{i\in \Z} V^i$ is a graded vector space over $\K$, then
$V[k]$ stands for the graded vector space with $V[k]^i:=V^{i+k}$. For $v\in V^i$ we set $|v|:=i$.
The phrase \emph{differential graded} is abbreviated by dg.
The $n$-fold symmetric product of a (dg) vector space $V$ is denoted  by $\odot^n V$, the full symmetric product space by  $\odot^\bullet V$.
For a finite group $G$ acting on a vector space $V$, we
denote via $V^G$ the space of invariants with respect to the action of $G$, and by $V_G$
the space of coinvariants $V_G = V/\{gv- v| v\in V, g\in G\}$. As we always work over a field $\K$ of characteristic zero, we have a canonical isomorphism $V_G\cong V^G$.

We  use freely the language of operads. For a background on operads we refer to the textbook \cite{LV}.
For an operad $\cP$ we denote by $\cP\{k\}$ the unique operad which has the following property:
for any graded vector space $V$ there is a one-to-one correspondence between representations of
$\cP\{k\}$ in $V$ and representations of
$\cP$ in $V[-k]$; in particular, $\cE nd_V\{k\}=\cE nd_{V[k]}$.

\bip

{\large
\section{\bf A variant of the Kontsevich graph complex}
}

\bip
\subsection{From operads to Lie algebras}
Let $\cP=\{\cP(n)\}_{n\geq 1}$ be an operad in the category of dg vector spaces with the partial compositions $\circ_i: \cP(n)\ot \cP(m) \rar \cP(m+n-1)$, $1\leq i\leq n$.
Then the map
$$
\Ba{rccc}
[\ ,\ ]:&  \sP \ot \sP & \lon & \sP\\
& (a\in \cP(n), b\in \cP(m)) & \lon &
[a, b]:= \sum_{i=1}^n a\circ_i b - (-1)^{|a||b|}\sum_{i=1}^m b\circ_i a
\Ea
$$
makes the vector space
$
\sP:= \prod_{n\geq 1}\cP(n)$  %%%\bigoplus_{n\geq 1}\cP(n)$
into a dg Lie algebra \cite{GV,KM}. Moreover, the Lie algebra structure descends to the subspace
of coinvariants $\sP_\bS:=  \prod_{n\geq 1}\cP(n)_{\bS_n}$. Via the identification of invariants and coinvariants $\sP_\bS \cong \sP^\bS$, we furthermore obtain a Lie algebra structure on the space of invariants
$\sP^\bS:=  \prod_{n\geq 1}\cP(n)^{\bS_n}$ as well.

% (via the composition $[\ ,\ ]:\wedge^2 \sP^\bS\rar \sP \rar \sP_\bS \cong \sP^\bS$) the Lie algebra structure on the subspace

\subsection{An operad of graphs and the Kontsevich graph complex}
For any integers $n\geq 1$ and $l\geq 0$ we denote by ${\sG}_{n,l}$ a set of graphs\footnote{ A {\em graph}\, $\Ga$ is, by definition,  a 1-dimensional $CW$-complex whose $0$-cells are called {\em vertices}\, and $1$-dimensional cells are called {\em edges}. The set of vertices of $\Ga$ is denoted by $V(\Ga)$ and the set of edges by $E(\Ga)$.}, $\{\Ga\}$, with $n$ vertices and $l$ edges
such that (i) the vertices of $\Ga$ are labelled by elements of $[n]:=\{1,\ldots, n\}$,
(ii) the set of edges, $E(\Ga)$, is totally ordered up to an even permutation.
For example, $\xy
(0,2)*{_{1}},
(7,2)*{_{2}},
 (0,0)*{\bullet}="a",
(7,0)*{\bu}="b",
\ar @{-} "a";"b" <0pt>
\endxy\in \sG_{2,1}$.
The group $\Z_2$ acts freely  on ${\sG}_{n,l}$ for $l\geq 2$ by changes of the total ordering;  its orbit
is denoted by $\{\Ga, \Ga_{opp}\}$. Let $\K\langle \sG_{n,l}\rangle$  be the vector space over  a field $\K$ spanned by isomorphism classes, $[\Ga]$, of elements of $\sG_{n,l}$ modulo the relation\footnote{Abusing notations we identify from now an equivalence class $[\Ga]$ with any
of its representative $\Ga$.}  $\Ga_{opp}=-\Ga$, and consider
a $\Z$-graded $\bS_n$-module,
$$
\sG \sr\sa (n):=\bigoplus_{l=0}^\infty \K\langle \sG_{n,l}\rangle[l].
$$
Note that graphs with two or more edges between any fixed pair of vertices do not contribute to
$\sG \sr\sa (n)$ so that we could have assumed right from the beginning that the sets $\sG_{n,l}$ do not contain graphs with multiple edges. The $\bS$-module, $\sG \sr\sa :=\{\sG \sr\sa (n)\}_{n\geq 1}$, is naturally an operad with the  operadic compositions given by
$$
\Ba{rccc}
\circ_i: & \sG \sr\sa (n)\ot \sG \sr\sa (m) & \lon &  \sG \sr\sa (m+n-1)\\
&  \Ga_1 \ot \Ga_2   &\lon & \sum_{\Ga\in \sG_{\Ga_1, \Ga_2}^i} (-1)^{\sigma_\Ga} \Ga
\Ea
$$
where $ \sG_{\Ga_1, \Ga_2}^i$ is the subset of $\sG_{n+m-1, \# E(\Ga_1) + \#E(\Ga_2)}$ consisting
of graphs, $\Ga$, satisfying the condition: the full subgraph of $\Ga$ spanned by the vertices labeled by
the set $\{i,i+1, \ldots, i+m-1\}$ is isomorphic to $\Ga_2$ and the quotient graph, $\Ga/\Ga_2$, obtained by contracting that subgraph to a single vertex, is isomorphic to $\Ga_1$. The sign $(-1)^{\sigma_\Gamma}$ is determined by the equality
$$
\bigwedge_{e\in E(\Ga)}e= (-1)^{\sigma_\Gamma}\bigwedge_{e'\in E(\Ga_1)}e' \wedge \bigwedge_{e''\in E(\Ga_2)}e''.
$$
The unique element in $\sG_{1,0}$ serves as the unit element in the operad  $\sG \sr\sa$. The associated  Lie algebra of  $\bS$-invariants, $((\sG \sr\sa\{-2\})^\bS,[\ ,\ ])$ is denoted, following notations
of \cite{Wi}, by $\ssf\sG\sC_2$. Its elements can be understood as
graphs from $\sG_{n,l}$ but with labeling of vertices forgotten, e.g.
$$
\xy
 (0,0)*{\bullet}="a",
(6,0)*{\bu}="b",
\ar @{-} "a";"b" <0pt>
\endxy = \frac{1}{2}\left(\xy
(0,2)*{_{1}},
(6,2)*{_{2}},
 (0,0)*{\bullet}="a",
(6,0)*{\bu}="b",
\ar @{-} "a";"b" <0pt>
\endxy + \xy
(0,2)*{_{2}},
(6,2)*{_{1}},
 (0,0)*{\bullet}="a",
(6,0)*{\bu}="b",
\ar @{-} "a";"b" <0pt>
\endxy\right)\in \ssf\sG\sC_2.
$$
The cohomological degree of a graph with $n$ vertices and $l$ edges is $2(n-1)-l$.
It is easy to check that $\xy
 (0,0)*{\bullet}="a",
(6,0)*{\bu}="b",
\ar @{-} "a";"b" <0pt>
\endxy$ is a Maurer-Cartan element in the Lie algebra $\ssf\sG\sC_2$. Hence we obtain a dg Lie algebra
$$
\left(\ssf\sG\sC_2, [\ ,\ ], d:=[\xy
 (0,0)*{\bullet}="a",
(6,0)*{\bu}="b",
\ar @{-} "a";"b" <0pt>
\endxy,\ ] \right).
$$
One may define a dg Lie subalgebra, $\sG\sC_2$, spanned by connected graphs with at least trivalent vertices and no edges beginning and ending at the same vertex. It is called
the {\em Kontsevich graph complex} \cite{Ko}.
We leave it to the reader to verify that the subspace $\sG\sC_2$ is indeed closed under both the differential and the Lie bracket.
We refer to \cite{Wi} for a detailed explanation of why studying
the dg Lie subalgebra  $\sG\sC_2$ rather than full Lie algebra $\ssf\sG\sC_2$ should be enough for most purposes. The cohomologies of $\sG\sC_2$ and $\ssf\sG\sC_2$ were partially computed in \cite{Wi}.

\subsubsection{\bf Theorem \cite{Wi}} \label{thm:negzero} (i) $
H^0(\sG\sC_2, d)\simeq \grt.$ (ii)  {\em For any negative integer}\, $i$, $H^i(\sG\sC_2, d)=0$.% while  is a one-dimensional vector space generated by the graph $
%\xy
%(0,-2)*{\bu}="A";
%(0,-2)*{\bu}="B";
%"A"; "B" **\crv{(6,6) & (-6,6)};
%\endxy
 %\in \ssf\sG\sC_2
%$.

%\mip

%It was realized by M. Kontsevich that the degree zero cocyles in $\sG\sC_2$
%This result implies that the Grothendieck-Teichm\"uller group $\GRT$ acts on the set of Poisson structures on an arbitrary manifold, up to homotopy.
%Note that to prove that action we need {\em both}\, results of the aforementioned Theorem.

\sip

We shall introduce  next a new graph complex which is responsible for the action of $\GRT$ on the set of quantum master functions on an odd symplectic supermanifold.

\subsection{A variant of the Kontsevich graph complex} The graph
$
\xy
(0,-2)*{\bu}="A";
(0,-2)*{\bu}="B";
"A"; "B" **\crv{(6,6) & (-6,6)};
\endxy
 \in \ssf\sG\sC_2
$
has degree $-1$ and satisfies
$$
[\xy
(0,-2)*{\bu}="A";
(0,-2)*{\bu}="B";
"A"; "B" **\crv{(6,6) & (-6,6)};
\endxy  ,   \xy
(0,-2)*{\bu}="A";
(0,-2)*{\bu}="B";
"A"; "B" **\crv{(6,6) & (-6,6)};
\endxy  ]=[\xy
(0,-2)*{\bu}="A";
(0,-2)*{\bu}="B";
"A"; "B" **\crv{(6,6) & (-6,6)};
\endxy , \xy
 (0,0)*{\bullet}="a",
(6,0)*{\bu}="b",
\ar @{-} "a";"b" <0pt>
\endxy]=0.
$$
Let $u$  be a formal variable of degree $2$ and consider the graph complex $\ssf\sG\sC_2[[u]]$ with the differential
$$
d_u:= d+ u \Delta,\ \ \ \ \ \mbox{where}\ \ \ \Delta:=[\xy
(0,-2)*{\bu}="A";
(0,-2)*{\bu}="B";
"A"; "B" **\crv{(6,6) & (-6,6)};
\endxy  , \ \ ].
$$
The subspace $\sG\sC_2[[u]]\subset \ssf\sG\sC_2[[u]]$ is a subcomplex of $(\ssf\sG\sC_2[[u]], d_u)$.

\begin{proposition}\label{2: prop on grt}
$H^0\left(\sG\sC_2[[u]], d_u\right)\simeq \grt$ and $H^{\leq -1}\left(\sG\sC_2[[u]]\right)=0$.
\end{proposition}
\begin{proof}
Consider a decreasing filtration of $\sG\sC_2[[\hu]]$ by the powers in $\hu$. The first term of the associated spectral sequence is
$$
\cE_1= \bigoplus_{i\in \Z} \cE_1^i,\ \ \ \ \cE_1^i=\prod_{p\geq 0} H^{i-2p}(\sG\sC_2, d) u^p
$$
with the differential equal to $u\Delta$.
As $H^0(\sG\sC_2, d)\simeq \grt$
and $H^{\leq -1}(\sG\sC_2, d)=0$, one gets the desired result.  %Similarly one gets the isomorphism
$H^0\left(\ssf\sG\sC_2[[u]], d_u\right)\simeq \grt$.

The projections $(\sG\sC_2[[u]],d_u)\to (\sG\sC_2,d)$ and $(\ssf\sG\sC_2[[u]],d_u)\to (\ssf\sG\sC_2,d)$ sending $u$ to 0 are maps of Lie algebras and induce isomorphisms in degree 0 cohomology. Since the isomorphisms of Theorem {\ref{thm:negzero}} (i) are maps of Lie algebras as shown in \cite{Wi}, so are the maps in the above Proposition.
\end{proof}

\subsection{Remark}\label{2:Remark on induction}  Let $\sigma$ be an element of $\grt$ and let $\Ga_\sigma^{(0)}$ be any cycle representing the cohomology class  $\sigma$  in the graph complex
$(\sG\sC_2, d)$.
Then one can construct a cocycle,
\Beq\label{cyclic_repr}
\Gamma^\hu_\sigma= \Ga_\sigma^{(0)}  + \Ga_\sigma^{(1)}\hu +  \Ga_\sigma^{(2)}\hu^2 +  \Ga_\sigma^{(3)}\hu^3+ \ldots,
\Eeq
 representing the cohomology class $\sigma\in\grt$ in the complex $\left(\sG\sC_2[[u]], d_u\right)$ by the following induction:
\sip

{\em 1st step}: As $d\Ga_\sigma^{(0)}=0$, we have $d (\Delta \Ga_\sigma^{(0)})=0$. As $H^{-1}( \sG\sC_2, d)=0$,
there exists $\Ga_\sigma^{(1)}$ of degree $-2$ such that $\Delta \Ga_\sigma^{(0)}= - d\Ga_\sigma^{(1)}$ and hence
$$
(d +u \Delta) \left(\Ga_\sigma^{(0)}  + \Ga_\sigma^{(1)}\hu\right)=0 \bmod O(\hu^2).
$$

{\em n-th step}: Assume we have constructed  a polynomial $\sum_{i=1}^n  \Ga_\sigma^{(i)}\hu^i$ such that
 $$
(d +\hu \Delta) \sum_{i=1}^n  \Ga_\sigma^{(i)}\hu^i=0 \bmod O(\hu^{n+1}).
$$
Then $d (\Delta \Ga_\sigma^{(n)})=0$, and, as  $H^{-2n-1}( \sG\sC_2, d)=0$, there exists a  graph
$\Ga_\sigma^{(n+1)}$ in $\sG\sC_2$ of degree $-2n-2$ such that $\Delta \Ga_\sigma^{(n)})=- d \Ga_\sigma^{(n+1)}$. Hence $(d +\hu \Delta) \sum_{i=1}^{n+1}  \Ga_\sigma^{(i)}\hu^i=0 \bmod O(\hu^{n+2})$.

%\sip
%Note that given the cocycle $\Gamma^\hu_\sigma$, then

%It is these $\hu$-dependent additions to $\Gamma_\sigma^{(0)}$ in $\Gamma^\hu_\sigma$  which make the action of $\GRT$ on  quantum master functions different from its  action \cite{Wi} on Poisson structures.

\bip

\bip

{\large
\section{\bf Quantum BV structures on odd symplectic manifolds}
}

% \bip
% \subsection{On $L_\infty$ automorphisms of Lie algebras} \label{3: subsec on Lie infty auto}

\bip

\subsection{Maurer-Cartan elements and gauge tranformations}
Let $(\fg=\oplus_{i\in \Z} \fg^i, [\ ,\ ], d)$ be a dg Lie algebra and consider the dg Lie algebra $\fg_\hbar := \hbar \fg[[\hbar]]=:\oplus_{i\in \Z} \fg^i_\hbar$, where $\hbar$ is a formal deformation parameter.
The group $G:=\exp(\fg_{\hbar}^0)$ (which is, as a set, $\fg_{\hbar}^0$ equipped with the standard Baker-Campbell-Hausdorff  multiplication) acts on $\fg^1_\hbar$,
$$
\ga \rar \exp(h)\cdot \ga:= e^{\mathrm{ad}_h}\ga -\frac{e^{\mathrm{ad}_h}-1}{\mathrm{ad}_h}dh,
$$
preserving its subset of Maurer-Cartan elements
$$
\cM\cC(\fg_\hbar) = \{\ga\in \fg^1_\hbar | d\ga + \frac{1}{2}[\ga,\ga]=0\}.
$$
We call the $G$-orbits in $\cM\cC(\fg_\hbar)$ the gauge equivalence classes of Maurer-Cartan elements.

 \sip

The group of $L_\infty$ automorphism of $\fg$ acts on $\cM\cC(\fg_\hbar)$ by the formula
$$
F\cdot \ga := \sum_{n\geq 1}\frac{1}{n!} F_n(\ga, \ldots, \ga)
$$
where $F_n$ is the $n$-th component of the $L_\infty$ morphism.
In particular, let $f$ be an $L_\infty$ derivation of $\fg$ without linear term. It exponentiates to an $L_\infty$ automorphism $\exp(f)$ of $\fg$, which acts on $\cM\cC(\fg_\hbar)$, and in particular on the set of gauge equivalence classes.
By a small calculation one may check that if we change $f$ by homotopy, i.~e., by adding $dh$ for some degree 0 element $h$ of the Chevalley-Eilenberg complex of $\fg$, then the induced actions of $\exp(f)$ and $\exp(f+dh)$ on the set of gauge equivalence classes agree.

\subsection{Quantum BV manifolds} Let $M$ be a $\Z$-graded manifold equipped with an odd symplectic structure $\om$ (of degree $1$). There always exist so called Darboux coordinates,
$(x^a, \psi_a)_{1\leq a\leq n}$, on $M$ such that $|\psi_a|=-|x^a| + 1$ and $\om=\sum_a dx^a\wedge d\psi_a$. The odd symplectic structure makes, in the obvious way, the structure sheaf  %$\f_M$
into a Lie algebra with brackets, $\{\ ,\ \}$, of degree $-1$. A less obvious fact is that $\om$ induces
a  degree $-1$ differential operator, $\Delta_\om$, on the invertible sheaf of semidensities, $\Ber(M)^{\frac{1}{2}}$ \cite{Kh}. Any choice of a Darboux coordinate system on $M$ defines an associated trivialization of the sheaf $\Ber(M)^{\frac{1}{2}}$; if one denotes the associated basis section
of  $\Ber(M)^{\frac{1}{2}}$ by $D_{x,\psi}$, then  any semidensity $D$ is of the form $f(x,\psi) D_{x,\psi}$ for some smooth function $f(x,\psi)$, and the operator $\Delta_\om$ is given by
 $$
 \Delta_\omega\left(f(x,\psi) D_{x,\psi}\right)=\sum_{a=1}^n %(-1)^{|x^a|}
  \frac{\p^2f}{\p x^a \p \psi_a} D_{x,\psi}.
 $$
Let $\hu$ be a formal parameter of degree $2$. A {\em quantum
master function}\ on $M$ is a $\hu$-dependent semidensity $D$ which satisfies the equation
$$
\Delta_\omega D=0
$$
and which admits, in some Darboux coordinate system, a form
$$
D=e^{\frac{S}{\hu}}D_{x,\psi},
$$
for some  $S\in \f_M[[\hu]]$ of total degree $2$, where $\f_M$ is the algebra of functions on $M$. In the literature it is this formal power series
in $\hu$ which is often called a quantum master function. Let us denote the set of all quantum master functions on $M$ by $\cQ\cM(M)$. It is easy to check that the equation $\Delta_\om D=0$ is equivalent to the following one,
\Beq\label{3:QME}
\hu\Delta S + \frac{1}{2}\{S,S\}=0,
\Eeq
where $\Delta:=\sum_{a=1}^n  \frac{\p^2}{\p x^a \p \psi_a}$. This equation is often called the {\em quantum master equation}, while a triple $(M, \om, S\in \cQ\cM(M))$ a {\em quantum BV manifold}.

\sip

Let us assume from now on that $M$ is affine or formal (i.~e., we work with $\infty$-jets of functions at some point) and that a particular Darboux coordinate system is fixed on $M$ up to affine transformations\footnote{This is not a serious loss of generality as any quantum master equation can be represented in the form (\ref{3:QME}). Our action of $\GRT$ on $\cQ\cM_\hbar(M)$
depends on the choice of an affine structure on $M$ in exactly the same way as
the classical Kontsevich's  formula for a universal formality  map \cite{Ko2} depends on such a choice. A choice of an appropriate affine connection on $M$ and methods of the paper \cite{D} can make our formulae for the $\GRT$ action invariant under the group of symplectomorphsims of $(M,\om)$; we do not address this {\em globalization}\, issue in the present note.} so that the algebra of function on $M$ is $\f_M\cong \K[x^a,\psi_a]$ or $\f_M\cong \K[[x^a,\psi_a]]$.

\sip

For later reference we will also consider solutions of \eqref{3:QME} that depend on a formal deformation parameter $\hbar$ of degree 0, $S\in \hbar \f_M[[u]][[\hbar]]$. We will call the set of such $S$ the \emph{set of formal solutions of the quantum master equation}\, and denote it by $\cQ\cM_\hbar(M)$.

\subsection{An action of $\GRT$ on quantum master functions}\label{3: subsection on GRT action on qmfunctions}
The constant odd symplectic structure on $M$ makes $\f_M$ into a representation
\Beq
\Ba{rccc}
\rho: & \sG \sr\sa(n) & \lon & \sE \mathsf n \mathsf d_V(n)=\Hom_{cont}(\f_M^{\ot n},\f_M)\\
      & \Ga &\lon & \Phi_\Ga
\Ea
\Eeq
of the operad $\sG \sr\sa$ as follows:
$$
\Phi_\Ga(S_1,\ldots, S_n) :=\pi\left(\prod_{e\in E(\Ga)} \Delta_e \left(S_1(x_{(1)}, \psi_{(1)},\hu)\ot S_2(x_{(2)}, \psi_{(2)},\hu)\ot \ldots\ot S_n(x_{(n)}, \psi_{(n)},\hu) \right)\right)
$$
where, for an edge $e$ connecting vertices labeled by integers $i$ and $j$,
$$
\Delta_e= \sum_{a=1}^n \frac{\p}{\p x_{(i)}^a}  \frac{\p}{\p \psi_{a(j)}}
+%  (-1)^{|x_{(j)}^a||\psi_{(i)a}|}
\frac{\p}{\p \psi_{a(i)}}  \frac{\p}{\p x_{(j)}^a}
$$
with the subscript $(i)$ or $(j)$ indicating that the derivative operator is to be applied to the $i$-th of $j$-th factor in the tensor product.
The symbol $\pi$ in \eqref{3: repr rho of Gra} denotes the multiplication map,
$$
\Ba{rccc}
\pi:&   V^{\ot n} & \lon & V\\
   & S_1\ot S_2\ot \ldots \ot S_n &\lon & S_1S_2\cdots S_n.
\Ea
$$

Let $V:=\f_M[[\hu]]$. Then by $\hu$-linear extension we obtain a continuous representation (in the category of topological $\K[[\hu]]$-modules)
\Beq \label{3: repr rho of Gra}
 \sG \sr\sa[[u]]  \lon  \sE \mathsf n \mathsf d_V=\Hom_{cont}(V^{\ot \cdot},V).
\Eeq

The space $V[1]$ is a topological dg Lie algebra with differential $\hu \Delta$ and
Lie bracket $\{\ ,\ \}$. These data define a Maurer-Cartan element, $\ga_{\cQ\cM}:=\hu\Delta \oplus
 \{\ ,\ \}$ in the Lie algebra $(\sE \mathsf n \mathsf d_V\{-2\})^\bS\subset CE^\bu(V,V)$, where $CE^\bu(V,V)$ is the Lie algebra of  coderivations
 $$
 CE^\bu(V,V)=\left(\mbox{Coder}(\odot^{\bu\geq 1}(V[2])), [\ ,\ ] \right)\ \ \mathrm{with}\ \ \ CE^\bu(V,V)_{(m)}:=\Hom(\odot^{\bu\geq m+1}(V[2]), V[2]),
 $$
of the standard graded co-commutative coalgebra, $\odot^{\bu\geq 1}(V[2])$, co-generated by a vector space $V$. The set $\cM\cC(CE^\bu(V,V))$ can be identified with the set of
 $L_\infty$ structures %\footnote{In our grading conventions the degree of $n$-th $L_\infty$ operation on $V$ is equal to $3-2n$.}
 on the space $V[1]$.
 %Any element $\ga\in \cM\cC(CE^\bu(V,V))$ defines a differential, $d_\ga:=[\ga,\ ]$ in $CE^\bu(V,V)$, and any element $g\in \Ker d_\ga\cap CE^0(V,V)_{(1)}$ defines an automorphism
% of  $\cM\cC(CE^\bu(V,V))$  which leaves the point $\ga$ invariant. Such an element $e^g$ can be interpreted as
% a $L_\infty$ automorphism of the $L_\infty$ algebra $(V[1],\ga)$. Moreover, if the cohomology class of $g$ in $H(CE^\bu(V,V), d_\ga)$ is non-trivial, then this $L_\infty$ automorphism is homotopy non-trivial.complex of .

\sip

The map sending an operad $\cP$ to the Lie algebra of invariants $\prod_n \cP\{-2\}(n)^{\bS_n}$ is functorial. Hence, from the representation \eqref{3: repr rho of Gra} we obtain a map of graded Lie algebras
$$
\ssf\sG\sC_2[[u]]\cong (\sG \sr\sa\{-2\}[[u]])^\bS \to (\sE \mathsf n \mathsf d_V\{-2\})^\bS \subset CE^\bu(V,V)
$$
One checks that the Maurer-Cartan element
$$ \xy
 (0,0)*{\bullet}="a",
(6,0)*{\bu}="b",
\ar @{-} "a";"b" <0pt>
\endxy
+ \hu \, \xy
(0,-2)*{\bu}="A";
(0,-2)*{\bu}="B";
"A"; "B" **\crv{(6,6) & (-6,6)};
\endxy \in \ssf\sG\sC_2[[\hu]]$$
is sent to the Maurer-Cartan element $\ga_{\cQ\cM}\in CE^\bu(V,V)$. Hence we obtain a morphism of dg Lie algebras
$$
 \left(\ssf\sG\sC_2[[\hu]], [\ ,\ ], d_h\right) \lon \left(CE^\bu(V,V), [\ ,\ ], \delta:=[\ga_{\cQ\cM},\ ]\right),
$$
and by restriction a morphism
$$
 \Phi\colon \left(\sG\sC_2[[\hu]], [\ ,\ ], d_h\right) \lon \left(CE^\bu(V,V), [\ ,\ ], \delta:=[\ga_{\cQ\cM},\ ]\right),
$$
Hence we also obtain a morphism of their cohomology groups,
$$
\grt\simeq H^0\left(\sG\sC_2[[\hu]], d_\hu\right) \lon H^0\left(CE^\bu(V,V), \delta\right).
$$
 Let $\sigma$ be an arbitrary element in $\grt$ and let
 $\Ga^\hu_\sigma$ be a cocycle representing $\sigma$ in the graph complex  $(\sG\sC_2[[\hu]], d_\hu)$.
 We may assume that $\Ga^\hu_\sigma$ consists of graphs with at least 4 vertices, see \cite{Wi}.
 Then the element $\Phi(\Ga_\sigma^\hu)$ describes an $L_\infty$ derivation of the Lie algebra $V[1]$ without linear term.
 By exponentiation we obtain an $L_\infty$ automorphism,
$$
F^\sigma=\left\{F^\sigma_n: \odot^n V \lon V[2-2n]\right\}_{n\geq 1},
$$
of the dg Lie algebra $(V[1], \hu\Delta, \{\ ,\ \})$  with $F^\sigma_1=\Id$.
Hence,  for any  formal quantum master function $S\in \cQ\cM_\hbar(M)$ the series
$$
S^\sigma:= S + \sum_{n\geq 2}\frac{1}{n!} F^\sigma_n(S, \ldots, S)
$$
gives again a formal quantum master function.\footnote{The series trivially converges since we work in the formal setting, i.~e., $S=\hbar(\cdots)$. Ideally, of course, one hopes to have a non-zero convergence radius in $\hbar$, but we cannot guarantee this.}
The induced action on gauge equivalence classes of such functions is well defined, i.~e., it does not depend on the representative $\Ga^\hu_\sigma$ chosen. %Furthermore, since the map $\grt\simeq H^0\left(\ssf\sG\sC_2[[\hu]], d_\hu\right)$ is a Lie algebra morphism, this action is indeed an action, i.~e., it respects the composition in $\GRT$.
This is the acclaimed homotopy action of $\GRT$ on  $\cQ\cM_\hbar(M)$
for any affine odd symplectic manifold $M$.

\begin{remark}
 As pointed out by one of the referees, there is also a stronger notion of ``homotopy action'' that holds in our setting. We will only consider the infinitesimal version.
 Then, we do not only have a Lie algebra morphism $\grt\to H^0\left(CE^\bu(V,V)\right)$, but an $L_\infty$ morphism $\grt\to CE^\bu(V,V)$ as follows. First, consider the truncated version $\left(\sG\sC_2[[\hu]]\right)^{tr}$ of the dg Lie algebra $\sG\sC_2[[\hu]]$, which is by definition the same as $\sG\sC_2[[\hu]]$ in negative degrees, zero in positive degrees, and consists of the degree zero cocycles in degree zero. By Proposition {\ref{2: prop on grt}} the canonical projection $\left(\sG\sC_2[[\hu]]\right)^{tr}\to \grt$ is a quasi-isomorphism. Hence we can obtain the desired $L_\infty$ morphism $\grt\to CE^\bu(V,V)$ by lifting the zig-zag
 \[
  \grt \stackrel{\sim}{\lon} \left(\sG\sC_2[[\hu]]\right)^{tr} \lon CE^\bu(V,V).
 \]
 This proves the first claim of the main Theorem.
\end{remark}

%It is perhaps, more natural in the present context to view
%\subsection{Remark} In QFT one often works with quantum master functions $S$ which are
%formal power series in the Darboux coordinates rather than real analytic functions. In that case
%one should view
%$u$ as a degree $0$ formal parameter so that the adjoint action of  ${e^{u\Phi_{\Ga_\sigma^\hu}}}\in G$ on $CE^\bu(V,V)[[u]]$ gives a continuous $L_\infty$ automorphism of the dg Lie algebra $(V[[u]][1], \hu\Delta, \{\ ,\ \})$, and hence induces a well-defined (in the sense, no problem with convergence) transformation of master functions from $\f_M[[\hu]][[u]]$.
%Such a transformation
%makes sense also for  $\Z_2$ graded (rather than $\Z$ graded)  quantum master functions; in that case
%one can simply set $u=\hu$.

\subsection{Remark}
It is a well known result due to D. Tamarkin \cite{T} that the Grothendieck Teichm\"uller group $\GRT$ acts on the operad of chains of the little disks operad.
In fact, one can show that this $\GRT$ action extends to an action on the operad of chains of the framed little disks operad, which is quasi-isomorphic to the Batalin-Vilkovisky operad. Hence one obtains in particular an action of $\GRT$ on the set of Batalin-Vilkovisky algebra structures on any vector space, and on their deformations, up to homotopy.
In our setting the algebra $\f_M$ is an algebra over the framed little disks operad. Any solution $S=S_0+u S_1+u^2 S_2+\cdots$ of the master equation \eqref{3:QME} yields a deformation of the Batalin-Vilkovisky structure on $\f_M$, up to homotopy. Concretely, to $S$ one may associate a $BV_\infty^{com}$-structure (see \cite{Kr} or \cite[section 5.3]{CMW}), whose $n$-th order ``BV'' operator is defined as $\Delta_n := [S_n,\cdot]$ (notation as in \cite[section 5.3]{CMW}).
The $\GRT$ action on solutions of the master equation described above can hence be seen as a shadow of this more general action of $\GRT$ on the framed little disks operad.
However, we leave the details to elsewhere.

\subsection*{Acknowledgements} We are grateful to K. Costello and to the anonymous referees for useful critical comments.

\def\cprime{$'$}

\end{document}